\documentclass[12pt,reqno]{amsart}

\usepackage{fullpage,amsfonts,amsmath,amssymb,amscd,youngtab}

\input xypic

\theoremstyle{plain}
\newtheorem*{lem}{Lemma}

\newtheorem*{cor}{Corollary}
\newtheorem*{thrm}{Theorem}
\newtheorem*{thm}{Theorem}

\theoremstyle{Definition}

\newtheorem*{defn}{Definition}

\newtheorem*{rem}{Remark}

\newcommand{\cO}{\mathcal{O}}
\newcommand{\irrep}{\textsf{Irrep}}

\newcommand{\Hom}{\mathrm{Hom}}

\newcommand{\e}{\varepsilon}

\newcommand{\C}{\mathbb{C}}

\renewcommand{\leq}{\leqslant}

\newcommand{\KZ}{\texttt{KZ}}

\DeclareMathOperator{\ext}{Ext}

\newcommand{\calH}{\mathcal{H}}
\newcommand{\cont}{\mathrm{cont}}

\newcommand{\h}{\mathfrak{h}}
\newcommand{\hreg}{\mathfrak{h}^{\mathrm{reg}}}

\begin{document}
\title{On category $\cO$ for the rational Cherednik algebra of $G(m,1,n)$: the
almost semisimple case}

\author{Richard Vale}

\address{Department of Mathematics, University
of Glasgow, Glasgow, G12 8QW, U.K.} \email{rv@maths.gla.ac.uk}

\date{\today}

\begin{abstract}
We determine the structure of category $\cO$ for the rational
Cherednik algebra of $G(m,1,n)$ in the case where the $\KZ$ functor satisfies a condition called \emph{separating simples}. As a consequence, we show that the property of having exactly $N-1$ simple modules, where $N$ is the number of simple modules of $G(m,1,n)$, determines the Ariki-Koike algebra up to isomorphism. 
\end{abstract}

\maketitle
\section{Introduction}
Let $W$ be the complex reflection group $G(m,1,n)=\mathbb{Z}_m \wr
S_n$ with its natural $n$--dimensional representation $\h$. It is
suggested in \cite{rouquiersurvey} that it may be possible to give a
complete description of the structure of category $\cO$ for the
rational Cherednik algebra of the group $W$, using generalisations
of the methods of \cite{Go}, \cite{BEG2}, \cite{GGOR} etc. for the
$m=1$ case. The aim of this paper is to do this for the case most
similar to \cite{BEG2}. Recall that in \cite{BEG2}, a complete
description of category $\cO$ was obtained in the case $G=S_n$ and
where the parameter $c$ belonged to $\frac{1}{h}+\mathbb{Z}_{\ge0}$.
In \cite{gmpn}, some of these results
were generalised to the case of $G(m,1,n)$ when the parameters are
chosen generically in a certain hyperplane. Here, we extend these
results, but giving instead a much cleaner condition involving the
$KZ$ functor, and then showing that when this condition holds, we
are essentially in the situation of \cite{gmpn}. Our main result is Theorem \ref{mainthm} below, which gives a complete description of category $\cO$ analogous to that proved in the $W=S_n$ case in \cite{BEG2}. As a corollary of Theorem \ref{mainthm}, we then prove that the Ariki-Koike algebra (ie. the Hecke algebra of $W$) is determined up to isomorphism by the property of having $|\irrep(W)|-1$ simple modules (see Corollary \ref{akassc}). Throughout, we
use the notation and definitions of \cite{gmpn}.
\section{Acknowledgements}
This work will form part of the author's PhD thesis at the University of Glasgow, funded by EPSRC. Part of this work was done while the author was visiting the University of Chicago, funded by the Leverhulme Trust. I would like to thank my supervisors K. A. Brown and I. Gordon, and also S. Ariki, A. Mathas and C. Stroppel for valuable discussions and for pointing out mistakes. 
\section{The rational Cherednik algebra}
\subsection{}In this section, we recall the basic facts about the
rational Cherednik algebra, before stating our main theorem. A
more detailed exposition can be found in \cite{gmpn}. We let
$W=G(m,1,n)$ with its reflection representation $\h$. Let $\{y_1,
\ldots, y_n\}$ be the standard basis of $\h$. With respect to this
basis, $W$ may be regarded as the group of $n \times n$ matrices
with exactly one nonzero entry in each row and column, the nonzero
entries being powers of $\e:=e^{\frac{2 \pi i}{m}}$. We also let
$\{x_1, \ldots, x_n\}$ denote the basis of $\h^*$ which is dual to
$\{y_1, \ldots, y_n\}$. \subsection{}The complex reflections in
$W$ are then the elements $s_i^t$ and $\sigma_{ij}^{(k)}$ defined
as follows: for $1 \le i \le n$ and $1 \le t \le m-1$, we define
$$s_i^t(y_i)=\e^ty_i, \quad s_i^t(y_j)=y_j, j \neq i$$
and for $1 \le i<j \le n$ and $0 \le k \le m-1$, define
$$\sigma_{ij}^{(k)}(y_i)=\e^{-k}y_j, \quad
\sigma_{ij}^{(k)}(y_j)=\e^{k}y_i, \quad
\sigma_{ij}^{(k)}(y_r)=y_r,  r \neq i,j.$$ Each of these elements
has a reflecting hyperplane $H$. The reflecting hyperplane of
$s_i^t$ is $\{y_i=0\}$ while the reflecting hyperplane of
$\sigma_{ij}^{(k)}$ is $\{y_i=\e^{-k}y_j\}$. Let $\mathcal{A}$ be
the set of these reflecting hyperplanes. For each $H \in
\mathcal{A}$, let $\alpha_H$ be a linear functional on $\h$ with
kernel $H$.
\subsection{} Let $\kappa=(\kappa_{00}, \kappa_0,\kappa_1, \ldots,
\kappa_{m-1})$ be a vector of complex numbers. Then the rational
Cherednik algebra $H_\kappa=H_\kappa(W)$ of $W$ is the quotient of
the $\C$--algebra $T(\h \oplus \h^*)* W$ by the relations
$[x_1,x_2]=0$ for $x_1,x_2 \in \h^*$, $[y_1,y_2]=0$ for $y_1,y_2
\in \h$, together with the commutation relations
\begin{multline*} [y,x] = y(x) +
\sum_{i=1}^ny(x_i)x(y_i)\sum_{j=0}^{m-1}(\kappa_{j+1}-\kappa_j)
\sum_{r=0}^{m-1}\e^{rj}s_i^r \\
+\kappa_{00}\sum_{1 \le i<j \le n}\sum_{k=0}^{m-1}y(x_i-\e^kx_j) x
(y_i-\e^{-k}y_j)\sigma_{ij}^{(k)}
\end{multline*}
for all $x \in \h^*$ and all $y \in \h$. In this paper, we will
take $\kappa_0=0$ throughout.
\subsection{The Dunkl representation}
Let $\hreg = \h \setminus (\cup_{H \in \mathcal{A}}H)$ and let
$\mathcal{D}(\hreg)$ denote the ring of differential operators on
$\hreg$. It is well--known (see for instance, \cite{DO},
\cite[Proposition 4.5]{EtGi}) that there is an injective
homomorphism
$$H_\kappa \hookrightarrow \mathcal{D}(\h^{\mathrm{reg}})*W$$
called the Dunkl representation. If $\delta = \prod_{H \in
\mathcal{A}} \alpha_H \in \C[\h]$, then $\C [\h^{\mathrm{reg}}] =
\C [\h]_\delta$ and the induced map
$$H_\kappa|_{\h^{\mathrm{reg}}}:=H_\kappa \otimes_{\C[\h]} \C [\h^{\mathrm{reg}}]
\rightarrow \mathcal{D}(\h^{\mathrm{reg}})*W$$ is an isomorphism
(\cite[Theorem 5.6]{GGOR}).
\subsection{Category $\cO$}
Following \cite{BEG1}, let $\mathcal{O}$ be the abelian category
of finitely-generated $H_\kappa$--modules $M$ such that
for $P \in \C[\h^*]^W$, the action of $P - P(0)$ is locally
nilpotent. Let $\textsf{Irrep}(W)$ denote the set of isoclasses of
simple $W$--modules. Given $\tau \in \irrep(W)$, define the
standard module $M(\tau)$ by:
$$M(\tau) = H_{\kappa} \otimes_{\C [\h^*]*W} \tau$$
where for $p \in \C[\h^*], w \in W$ and $v \in \tau$, $pw \cdot v
:= p(0) wv$.

 In \cite{DO}, it is proved that $M(\tau)$ has a
unique simple quotient $L(\tau)$, and \cite{GGOR} prove that $\{
L(\tau) | \tau \in \irrep(W) \}$ is a complete set of
nonisomorphic simple objects of $\cO$, and that every object of
$\cO$ has finite length. Furthermore, it is proved in \cite{GGOR} that category $\cO$ is a \emph{highest weight category} in the sense of \cite{CPS}. In particular, every simple object $L(\tau)$ of $\cO$ has a projective cover $P(\tau)$ and an injective envelope $I(\tau)$, and \emph{BGG reciprocity} holds, that is, $[P(\tau):L(\sigma)]=[M(\sigma):L(\tau)]$ for all $\sigma, \tau$.
\subsection{The KZ functor}
The group $B_W:= \pi_1 (\h^{\mathrm{reg}}/W)$ is called the braid
group of $W$. In \cite{GGOR}, a functor
$$\texttt{KZ} : \cO \rightarrow \C B_W-\mathrm{mod}$$
is constructed as follows. If $M \in \cO$ then
$M|_{\h^{\mathrm{reg}}}:= \C [\h^{\mathrm{reg}}] \otimes_{\C[\h]} M$
is a finitely-generated module over $\C [\h^{\mathrm{reg}}]
\otimes_{\C[\h]} H_\kappa \cong \mathcal{D}(\h^{\mathrm{reg}})*W$.
In particular, $M$ is a $W$--equivariant $\mathcal{D}$--module on
$\h^{\mathrm{reg}}$ and hence corresponds to a $W$--equivariant
vector bundle on $\h^{\mathrm{reg}}$ with a flat connection
$\nabla$. The horizontal sections of $\nabla$ define a system of
differential equations on $\h^{\mathrm{reg}}$ which, by a process
described in \cite{rouquiersurvey}, give a monodromy
representation of $\pi_1 (\h^{\mathrm{reg}}/W)$. By definition,
$\texttt{KZ}(M)$ is the monodromy representation of $\pi_1
(\h^{\mathrm{reg}}/W)$ associated to $M$.
\subsection{}
\label{kznumbers} By \cite[Section 5.25]{GGOR},
the monodromy representation factors through the Hecke algebra
$\mathcal{H}$ of $W$. This is the quotient of $\C B_W$ by relations
given in \cite[Section 5.2,5]{GGOR}. From the braid diagram in
\cite[Table 1]{BMR}, we see that $\mathcal{H}$ is generated by
$T_{s}, T_{t_2}, \ldots , T_{t_m}$ subject to the relations:
{\allowdisplaybreaks \begin{align*}
T_{s}T_{t_2}T_sT_{t_2} - T_{t_2}T_s T_{t_2} T_s &=0 \\
[T_s, T_{t_i} ] &= 0 & i \ge 3 \\
T_{t_i} T_{t_{i+1}} T_{t_i} - T_{t_{i+1}} T_{t_i}
T_{t_{i+1}} &=0 & 2 \le i \le r \\
[T_{t_i}, T_{t_j} ] &= 0 & |i-j|> 1 \\
(T_{t_i} -1 )(T_{t_i} + e^{2 \pi i \kappa_{00} }) &= 0 & 2 \le i \le r \\
(T_s -1 ) \prod_{j=1}^{m-1} (T_s - \e^{-j} e^{-2 \pi i \kappa_j} ) &= 0  \\
\end{align*}}
We see that $\mathcal{H}$ is the Ariki-Koike algebra of \cite{AK}
with parameters $q= e^{2 \pi i \kappa_{00}}$, and $u_i =
\e^{-(m-i+1)} e^{-2 \pi i \kappa_{m-i+1}}$ for $1 \le i \le m$,
where as before, $\e = e^{2 \pi i/m}$. Note in particular that
$u_i \neq 0$ for all $i$.
\subsection{}\label{otor}
Therefore, $\texttt{KZ}$ gives a functor $\texttt{KZ}: \cO
\rightarrow \mathcal{H}-\mathrm{mod}$. By \cite[Section
5.3]{GGOR}, $\texttt{KZ}$ is exact, and if $\cO_{\mathrm{tor}}$ is
the full subcategory of those $M$ in $\cO$ such that
$M|_{\h^{\mathrm{reg}}} =0$ then $\texttt{KZ}$ gives an
equivalence $\cO/\cO_{\mathrm{tor}} \tilde{\rightarrow}
\mathcal{H}-\mathrm{mod}$ \cite[Theorem 5.14]{GGOR}.
\section{A condition on KZ}
\subsection{}
Our aim is to study category $\cO$ in the situation where it is,
in some sense, as close as possible to being semisimple. We make
the following definition:
\begin{defn}
Say $\KZ : \cO \rightarrow \mathcal{H}-\mathrm{mod}$
\emph{separates simples} if whenever $S \ncong T$ are simple
objects of $\cO$, then $\KZ(S) \ncong \KZ(T)$.
\end{defn}
\subsection{}
Now we state the main theorem.
\begin{thrm}\label{mainthm}
Suppose $m>1$ and $n>1$ and $\KZ$ separates simples. Then either
$\cO$ is semisimple, or the following hold:
\begin{enumerate}
\item There exists a linear character $\chi$ of $W$ such that
$L(\chi)$ is finite-dimensional and all the other simple objects in
$\cO$ are infinite-dimensional. \item There exists a positive
integer $r$ not divisible by $m$, such that
$\mathrm{dim}L(\chi)=r^n$.
\item Let $q \in \mathbb{N}$ be the residue of $r$ modulo $m$, $1 \le q \le m-1$. Then
there is a representation $\h_q$ of $W$ with $\dim\h_q = \dim\h$
such that if $\tau \notin \{\wedge^i \h_q \otimes \chi | \: 0 \le i
\le n \}$ , then $M(\tau)=L(\tau)$. \item $\cO = \cO^\wedge \oplus
\cO^{ss}$ where $\cO^\wedge$ is generated by the $L(\wedge^i \h_q
\otimes \chi)$ and $\cO^{ss}$ is a semisimple category generated by
the other simple objects. \item The composition multiplicities in
$\cO^\wedge$ are $$[M(\wedge^i \h_q \otimes \chi) : L(\wedge^j \h_q
\otimes \chi)] = \begin{cases}  1
& \text{if $j=i, i+1$} \\
 0  & \text{otherwise}  \end{cases}$$
\end{enumerate}
\end{thrm}
\subsection{}
Before proving Theorem \ref{mainthm}, we make some remarks. Theorem
\ref{mainthm} is an analogue for $G(m,1,n)$ of various results of
the papers \cite{BEG2} and \cite{Go}. In fact, in \cite{BEG2} it is
shown that whenever $H_\kappa(S_n)$ has a finite-dimensional module,
then all but one of the simple modules in category $\cO$ are
infinite-dimensional and the structure of category $\cO$ is similar
to the result of Theorem \ref{mainthm} (see \cite[Theorem 1.2,
Theorem 1.3]{BEG2}). Although the methods we use for proving Theorem
\ref{mainthm} are similar to those of \cite{BEG2}, we have to use
different arguments to get round the problem that in the $G(m,1,n)$
case, the functor $\texttt{KZ}$ is not known to take standard
modules $M(\lambda)$ in $\cO$ to the corresponding Specht modules
$S^\lambda$ for $\mathcal{H}$, even on the level of Grothendieck
groups. We also have to do some work to calculate the blocks of the
Hecke algebra at the parameters that we are interested in.
\subsection{}
One reason why Theorem \ref{mainthm} is of interest is that it
gives a source of examples of choices of $\kappa$ such that there
is a finite-dimensional object in category $\cO$, and yet category
$\cO$ is completely understood.
\subsection{}
The rest of this section is devoted to the proof of Theorem
\ref{mainthm}. The proof proceeds as follows. In Section \ref{ak},
we recall some facts about the representations of the Ariki-Koike
algebra. We use these facts in Section \ref{beginning} to Section
\ref{parts1and2} to prove parts (1) and (2) of Theorem
\ref{mainthm}. Next, between Section \ref{blocks} and Section
\ref{parts3and4}, we compute the blocks of the Ariki-Koike algebra
in our situation by a combinatorial argument. This enables us to
prove parts (3) and (4) of Theorem \ref{mainthm}. Finally, in
Sections \ref{compo} and \ref{part5}, we prove part (5) of Theorem
\ref{mainthm}.
\section{The Ariki-Koike algebra}\label{ak}
\subsection{} Let us recall some facts about the Ariki-Koike
algebra. This is the algebra $\mathcal{H}$ introduced in Section
\ref{kznumbers}, also called the Hecke algebra of $W$. It depends
on parameters $q, u_1, \ldots, u_m \in \C$ and we are only
interested in the case where these parameters are all nonzero.
\subsection{} We use the following conventions. For us, a
partition of $n$ is a sequence $\lambda_1 \ge \lambda_2 \ge \cdots
\ge \lambda_k$ with $\sum \lambda_k =n$. A partition $\lambda$ will
be identified with its Young diagram, and we use the non-Francophone
convention for Young diagrams. That is, the Young diagram of
$\lambda$ has $\lambda_i$ boxes in row $i$, row 1 being the top row.
A multiparition of $n$ is an $m$--tuple $(\lambda^{(1)}, \ldots,
\lambda^{(m)})$ where the $\lambda^{(i)}$ are partitions with $\sum
|\lambda^{(i)}| =n$. Following the paper \cite{arikimathas}, we may
regard a multiparition as a subset of $\mathbb{N} \times \mathbb{N}
\times \mathbb{N} $ by thinking of it as an $m$--tuple of Young
diagrams. A \emph{node} is any box of $\lambda$. More generally, a
node will be any element of $\mathbb{N} \times \mathbb{N} \times
\mathbb{N} $.
\subsection{} It has been shown (see \cite{mathasbook}) that for
each multipartition $\lambda = (\lambda^{(1)}, \ldots,
\lambda^{(m)})$ of $n$, there is a \emph{Specht module}
$S^\lambda$ for $\mathcal{H}$. Each $S^\lambda$ has a quotient
$D^\lambda$ which is either 0 or simple. The set
$\{D^\lambda|D^\lambda \neq 0\}$ is a complete set of
nonisomorphic simple $\mathcal{H}$--modules. We will need a
parametrisation of this set. There are 2 different
parametrisations, depending on whether $q=1$ or $q \neq 1$.
\subsection{}\label{q=1case}
 If $q=1$ then \cite[Theorem 3.7]{mathasaksimples}
states that $D^\lambda \neq 0$ if and only if $\lambda^{(s)}=0$
whenever $s<t$ and $u_s=u_t$.
\subsection{}
If $q \neq 1$ then the description, due to Ariki and
stated in \cite[Theorem 3.24]{mathassurvey} is more complicated.
The nonzero $D^\lambda$ are in bijection with the set of
\emph{Kleshchev multipartitions}, which we now describe.

Given a multipartition $\lambda$, the \emph{residue} of a node $x$
in row $i$ and column $j$ of $\lambda^{(k)}$ is defined to be
$u_kq^{j-i}$. A node $x$ in $\lambda$ with residue $a$ is called a
\emph{removable $a$--node} if $\lambda \setminus \{ x\}$ is a
multipartition. A node $x$ not in $\lambda$ with residue $a$ is
called an addable $a$--node if $\lambda \cup \{ x\}$ is a
multipartition.

Say a node $y \in \lambda^{(\ell)}$ is \emph{below} a node $x \in
\lambda^{(k)}$ if either $\ell >k$, or $\ell=k$ and $y$ is in a
lower row than $x$.

A removable $a$--node $x$ is called \emph{normal} if whenever
${x'}$ is an addable $a$--node below $x$ then there are more
removable $a$--nodes between $x$ and ${x'}$ than there are
addable $a$--nodes. The highest normal $a$ node in $\lambda$ is
called the \emph{good} $a$--node.

The set of Kleshchev multiparitions is defined inductively as
follows: $\O$ is Kleshchev, and otherwise $\lambda$ is Kleshchev
if and only if there is some $a \in \C$ and a good $a$--node $x
\in \lambda$ such that $\lambda \setminus \{ x\}$ is Kleshchev.
\subsection{}\label{ihategrojnowski}
Finally we need a description of the blocks of $\mathcal{H}$. This
is given in \cite[Corollary 2.16]{lylemathas}. Recall that the Specht
modules are partitioned into blocks as follows: two Specht modules
$S^\lambda$ and $S^\mu$ are in the same block if and only if there
is a sequence $S^{\lambda_1}, S^{\lambda_2}, \ldots, S^{\lambda_t}$
with $S^{\lambda_1}=S^\lambda$, $S^{\lambda_t}=S^\mu$ and such that
$S^{\lambda_i}$ and $S^{\lambda_{i+1}}$ have a common composition
factor for all $i$. Define the \emph{content} $\cont(\lambda)$ of a
multipartition $\lambda$ to be the multiset of residues of
$\lambda$, ie. the set of residues counted according to
multiplicity. Then for $q \neq 1$, two Specht modules $S^\lambda$
and $S^\mu$ are in the same block if and only if
$\cont(\lambda)=\cont(\mu)$.

\section{Proof of Theorem \ref{mainthm}}
\subsection{}\label{beginning}
To begin the proof, suppose $\KZ$ separates simples. If $\cO$ is
not semisimple, then $\mathcal{H}$ is not semisimple, and we claim there
exists $S \in \cO$ with $\KZ(S)=0$. Indeed, if $\KZ(S) \neq 0$ for all simple objects $S \in \cO$ then $\calH$ has $|\irrep(W)|$ simple modules, but it is well-known that this implies that $\calH$ is semisimple. We give here a proof using the Cherednik algebra.
\begin{lem}
Suppose $\calH$ has $|\irrep(W)|$ simple modules. Then $\calH$ is semisimple.
\end{lem}
\begin{proof}
Since $\KZ$ is exact, each simple $\calH$--module is the image of a simple object of $\cO$ under $\KZ$. Therefore, the category $\cO_{\mathrm{tor}} \subset \cO$ is $0$. So $\KZ$ induces an equivalence $\cO \rightarrow \calH-\mathrm{mod}$. We show that $\cO$ is a semisimple category. By \cite[(32)]{DO}, there is an ordering $\leq$ on $\irrep(W)$ such that $[M(\tau):L(\sigma)] \neq 0$ implies $\tau \leq \sigma$. By \cite[Proposition 5.2.1]{GGOR}, if $L(\sigma)|_{\hreg} \neq 0$ then $L(\sigma) \subset M(\tau)$ for some $\tau$. Combining this fact with induction on the ordering $\leq$ yields $M(\tau)=L(\tau)$ for all $\tau$. But it is observed in \cite[Remark following Lemma 2.12]{BEG1} that $M(\tau)=L(\tau)$ for all $\tau$ if and only if $\cO$ is semisimple. Since there is an equivalence of categories $\cO \cong \mathcal{H}-\mathrm{mod}$, $\mathcal{H}-\mathrm{mod}$ is a semisimple category and so $\calH$ is a semisimple algebra. 
\end{proof}
\begin{rem}
Note that the above proof works for any complex reflection group $W$, where $\calH$ is the Hecke algebra of $W$ as defined in \cite[Section 5.2.5]{GGOR}.
\end{rem}
Therefore we have shown that if $\cO$ is not semisimple then there is some simple $S \in \cO$ with $\KZ(S)=0$, and $\KZ(T) \neq 0$ for all simples $T \ncong S$ by our assumption on
$\KZ$. Since $\KZ$ separates simples, we also have that $\#\{ \KZ(T) : T \: \text{simple}, T \ncong S\}= |\irrep(W)|-1$. Furthermore, if $T$ is simple then so is $\KZ(T)$, because $\KZ$ induces an equivalence
$\cO/\cO_{\mathrm{tor}} \rightarrow \mathcal{H}-\mathrm{mod}$, and
the localisation to $\h^\mathrm{reg}$ preserves simple objects. So
$\mathcal{H}$ has exactly $|\textsf{Irrep}(W)|-1$ simple modules.

\subsection{}
Next, we show that $q \neq 1$. Suppose $q=1$. Then by Section
\ref{q=1case}, since $\mathcal{H}$ is not semisimple, there must
be some $s<t$ with $u_s=u_t$. Under the assumption that $n>1$,
there are at least 3 multipartitions $\lambda$ with $\lambda^{(s)}
\neq \O$. Hence, there are at least 3 $D^\lambda$ which are zero
and so $\mathcal{H}$ cannot have $|\textsf{Irrep}(W)|-1$ simple
modules. So $q \neq 1$. Therefore, the simple
$\mathcal{H}$--modules are in bijection with Kleshchev
multipartitions.

\subsection{}
Ariki's semisimplicity criterion \cite{arikisemisimple} tells us
that $[n]_q!\prod_{i < j}\prod_{-n < c < n} (u_i - q^c u_j) =0$.
Therefore, either there are $i,j,c$ with $u_i = q^cu_j$, or else
$[n]_q!=0$. We show $[n]_q! \neq 0$. Suppose that $[n]_q!=0$. Then
there is a $k$, $1 \le k \le n$ with $q^k=1$ and $q^\ell \neq 1$,
$0 < \ell < k$. Since $q \neq 1$, the simple
$\mathcal{H}$--modules are in bijection with Kleshchev
multipartitions. Let $\rho_k$ be the partition of $k$ with one
part, ie. the Young diagram of $\rho_k$ is a row of $k$ boxes.
Then $\rho_k$ is not Kleshchev, because the only removable node of
$\rho_k$, call it $\mu$, cannot be good, because it is not normal.
Indeed, the node labelled $\lambda$ in the diagram below is an
addable node below $\mu$ with the same residue as $\mu$, and there
are no removable nodes between them.
$$\stackrel{\text{$k$ boxes in row}}{\overbrace{\young(\hfil \hfil
\hfil \hfil \mu,\lambda)}}.$$ Hence, $\rho_k$ is not Kleshchev and
therefore $\rho_n$, a row of $n$ boxes, is not Kleshchev. So we
may define multipartitions $\lambda_1 = (\rho_n , \O , \ldots,
\O)$ and $\lambda_2 = (\O, \rho_n , \O , \ldots, \O)$, neither of
which is Kleshchev (here we use the hypothesis that $m>1$). This
contradicts the fact that there is only one non-Kleshchev
multipartition, and so $[n]_q! \neq 0$.

\subsection{}
Therefore, there exist integers $1 \le i,j \le n$ and $-n < c <n$
such that $u_i = q^c u_j$. Writing what this means in terms of the
$\kappa_i$, we get
$$m(\kappa_j -\kappa_i) -mc\kappa_{00} - (i-j) \in m \mathbb{Z}.
$$
The next step is to show that $|c|=n-1$.

\subsection{}
Redefining $c$ if necessary, we have that there are $i < j$ with
$q^c u_i = u_j$. Either $c \ge 0$ or $c \le 0$. Consider the case
$c \ge 0$. In this case, let $\rho_{c+1}$ be a row of $c+1$ boxes,
and take a multipartition $\tau$ with $\rho_{c+1}$ as its
$i^\mathrm{th}$ part and $\o$ everywhere else. If $c<n-1$ then
consider two multipartitions defined as follows: $\lambda$ is the
multipartition of $n$ whose $i^\mathrm{th}$ part is $\rho_n$ and
$\mu$ is the multipartition of $n$ whose $i^\mathrm{th}$ part is
$$\stackrel{\text{$n-1$ boxes in row}}{\overbrace{\yng(5,1)}}.$$
Then $\tau$ is not Kleshchev, and so $\lambda$ is clearly not
Kleshchev. Also, $\mu$ is not Kleshchev, essentially because
 $\mu \supset \tau$ (note that, even after some nodes have been removed
 from $\mu$, the node at the right hand end of $\tau$
 can never be a good node, since we have established that $q^{c+1} \neq 1$).
  Hence there are 2
non-Kleshchev multipartitions, which contradicts our hypothesis.
So $c=n-1$.

In the $c \le 0$ case, we take $\gamma_{c+1}$ to be a column of
$-c+1$ boxes, and do a similar argument to show that $c=-(n-1)$.

\subsection{}\label{qorder}
The above argument shows that the mulitplicative
order of $q$ must be at least $2n-1$. Indeed, suppose $q^{n+a}=1$
where $a$ is a nonnegative integer. Then if $q^{n-1}u_i = u_j$ for
some $i, j$, we get $q^{-a-1}u_i = u_j$. But the above argument in
the $c \le 0$ case shows that $-a-1 \le -n$ or else we would have
more than one non-Kleshchev multipartition.

\subsection{}\label{parts1and2}
Now we may rewrite our condition on the parameters as
$$m(\kappa_j-\kappa_i)+(-1)^am(n-1)\kappa_{00} = (i-j)+mt$$
for some $a \in \{0,1\}$ and some $t \in \mathbb{Z}$. Note that
$(i-j) +mt$ cannot be zero because $1 \le i,j \le m$. If it is
positive, multiply through by $-1$ (possibly interchanging the
roles of $i$ and $j$, and changing $a$), so assume that $(i-j)+mt
<0$. Now we do a so-called twist. Consider the linear character of
$W$ which sends $\sigma_{ij}^{(\ell)}$ to $(-1)^a
\sigma_{ij}^{(\ell)}$ for all $i,j, \ell$, and which sends $s_k$
to $\e^{-j}s_k$. Explicitly checking with a set of generators and
relations of $W$ shows that this is a well-defined character of
$W$. Now by \cite[Section 5.4.1]{GGOR}, we have an isomorphism of
Cherednik algebras $\psi: H_\kappa \rightarrow H_{\kappa'}$
where $\kappa_{00}' =(-1)^a \kappa_{00}$ and $\kappa_{j}' =
\kappa_{j+i}-\kappa_i$. (These equations for $\kappa_i'$ follow
from writing down the generators and relations for
$H_{\kappa'}$). The twist $\psi$ induces an auotequivalence of
category $\cO$ which preserves the dimension of the objects
\cite[Section 5.4.1]{GGOR}. Our new parameters satisfy
$$m \kappa_{j-i}' + m(n-1)\kappa_{00}' = (i-j) + mt <0.$$
Now we are in a position where we can use \cite[Section 4.1]{CE}.
Translating our parameters into the language of \cite{CE}, we get
$$m (n-1) k + 2\sum_{j=1}^{m-1}c_j \frac{1-\e^{-jq}}{1-\e^{-j}}=r$$
where $r = (j-i)-mt$ is a positive integer of the form $(p-1)m +q$
for some nonnegative integer $p$ and some $1 \le q \le m-1$ (of
course, this $q$ is not the same $q$ which was a parameter in the
Hecke algebra). Then we have the module $\tilde{Y_c}$ which is a
quotient of $M(\mathsf{triv})$. Furthermore, since $[n]_q! \neq
0$, we may apply \cite[Theorem 4.3]{CE} to conclude that
$\tilde{Y_c}$ is finite-dimensional. Therefore, $L(\mathsf{triv})$
is finite-dimensional. By \cite[Section 5.4.1]{GGOR}, twisting by
$\psi$ sends $L(\chi)$ to $L(\mathsf{triv})$ for some linear
character $\chi$ of $W$. Furthermore, $\dim L(\chi) = \dim
L(\mathsf{triv}) = r^n$ by \cite[Theorem 2.3 (iii)]{CE}. Since
$L(\chi)$ is finite-dimensional, $\KZ(L(\chi))=0$, and therefore
$\KZ(L(\tau)) \neq 0$ for $\tau \neq \chi$, by our assumption that
$\KZ$ separates simples. Therefore $L(\tau)$ is
infinite-dimensional if $\tau \neq \chi$. We have proved parts
$(1)$ and $(2)$ of Theorem \ref{mainthm}.

\subsection{Blocks}\label{blocks}
To proceed further, it  is necessary to calculate the blocks of
the Hecke algebra.
\subsection{\emph{Standing assumption}}\label{stand}
We have parameters $q$ and $u_1, \ldots, u_m$ for the Hecke
algebra. We are assuming that there is exactly $1$ non-Kleshchev
multipartition, and we have already shown that $q^{n-1}u_i =u_j$
for some $i \neq j$. We have shown that under this condition on
the parameters, the unique non-Kleshchev multipartition has a row
of $n$ boxes as its $i^\mathrm{th}$ part, and all other parts
$\O$. First, we prove the following lemma.
\begin{lem}\label{lem1}
If $k \neq i,j$ then for each $\ell \neq k$, we have $u_k/u_\ell
\neq q^c$ for any $-n<c<n$.
\end{lem}
\begin{proof}
Suppose $u_k=q^cu_{\ell}$. If $\ell \neq i,j$ then it would follow
from the earlier calculations that there is another non-Kleshchev
multipartiton, so we need only consider the case where $\ell=i$ or
$\ell=j$. Suppose $i<j$. If $\ell=i$ then suppose there is
$-n<c<n$ with $u_k=q^cu_i$, and $u_j=q^{n-1}u_i$. If $c<0$ then
considering a multipartition whose only nontrivial part is a
column $\gamma_n$ or a row $\rho_n$ in the $i^\mathrm{th}$
position, we have that there is more than one non-Kleshchev
multipartition. On the other hand, if $c \ge 0$ then $u_k = q^c
u_i= q^{c-(n-1)}u_j$ and hence there exists a non-Kleshchev
multipartition which is $\O$ except in the $j^\mathrm{th}$
position, and one which is $\O$ except in the $i^\mathrm{th}$
position. Similarly, if $\ell=j$, we reach the same conclusion,
and so such a $c$ cannot exist. Similar arguments deal with the
$i>j$ case.
\end{proof}

\subsection{}
Recall from Section \ref{ihategrojnowski} that if $\alpha$ and
$\beta$ are multipartitions then the Specht modules $S^\alpha$ and
$S^\beta$ belong to the same block if and only if
$\cont(\alpha)=\cont(\beta)$. The next lemma is needed to study
the content of a multipartition.
\begin{lem}\label{lem2}
Under the assumptions of Section \ref{stand}, let
$\alpha=(\alpha^{(1)}, \alpha^{(2)}, \ldots, \alpha^{(m)})$ be a
multipartition of $n$. Then $\cont(\alpha^{(r)})\cap
\cont(\alpha^{(s)})=\O$ for all $r \neq s$.
\end{lem}
\begin{proof}
By Lemma \ref{lem1} and our assumption that $q^{n-1}u_i=u_j$, we
get that for all $r,s$, $u_r/u_s \neq q^c$ for any $-(n-1)<c<n-1$.
Now, if the residue of some node $x$ in $\alpha^{(r)}$ is equal to
the residue of some other node $y$ in $\alpha^{(s)}$, then
$$u_rq^{\mathrm{col}(x)-\mathrm{row}(x)}=u_sq^{\mathrm{col}(y)-
\mathrm{row}(y)}.$$ But if $t:=
\mathrm{col}(x)+\mathrm{row}(y)-\mathrm{row}(x)-\mathrm{col}(y)$
then $u_s/u_r = q^t$ but $t \leq n-2$ and $t \ge -(n-2)$, a
contradiction.
\end{proof}
\subsection{}
The next lemma is useful in determining a multipartition from its
content.
\begin{lem}\label{lem3}
Under the assumptions of Section \ref{stand}, if $\alpha$ and
$\beta$ are multipartitions of $n$ and $1 \le k \le m$, then
$\cont(\alpha^{(k)})=\cont(\beta^{(k)})$ implies
$\alpha^{(k)}=\beta^{(k)}$.
\end{lem}
\begin{proof}
We show that if two nodes of $\alpha^{(k)}$ have the same residue,
then they lie on the same diagonal. It will follow that the
multiplicity of a residue in $\cont(\alpha)$ is equal to the
length of the corresponding diagonal of $\alpha$. The same is true
of $\beta$. Thus under the hypothesis, the Young diagrams $\alpha$
and $\beta$ have diagonals of the same lengths, thus they are
equal.

Suppose then that nodes $(i,j)$ and $(i^{'},j^{'})$ in
$\alpha^{(k)}$ have the same residue. Then
$u_kq^{j-i}=u_kq^{j^{'}-i^{'}}$. Thus $q^{j-i-j^{'}+i^{'}}=1$ and
therefore if $j-i \neq j^{'}-i^{'}$ then either
$z:=j-i-j^{'}+i^{'} \ge n$ or $z \le -n$. But $j+i^{'}, j^{'}+i
\le n+1$ and so $z$ cannot be either greater than $n$ or less than
$-n$. Therefore, $z=0$ and $j-i=j^{'}-i^{'}$. In other words,
$(i,j)$ and $(i^{'},j^{'})$ lie on the same diagonal.
\end{proof}
%
%
%
%
\subsection{}\label{blocksclaim}
We are finally in a position to calculate the blocks of the Hecke
algebra. In order to determine the blocks of $\mathcal{H}$, we
first note that if $\rho_a$ denotes a row of length $a$ and
$\gamma_b$ a column of length $b$, then we may define a
multipartition $\lambda_{a}$ to have $\rho_a$ in the
$i^\mathrm{th}$ place and $\gamma_{n-a}$ in the $j^\mathrm{th}$
place. For example, if $m=3$, $n=3$, $i=3$, $j=2$ then
$${\tiny\Yvcentermath1{\lambda_0=\left(
\begin{array}{ccc} \O & \yng(1,1,1) & \O \end{array} \right),
\lambda_1=\left(\begin{array}{ccc} \O & \yng(1,1) & {\yng(1)}
\end{array} \right) , \lambda_2=\left(\begin{array}{ccc} \O &
{\yng(1)} & {\yng(2)} \end{array} \right), \lambda_3=
\left(\begin{array}{ccc} \O & \O & {\yng(3)} \end{array}
\right)}}.$$ Then if $q^{n-1}u_i =u_j$, then
$\mathrm{cont}(\lambda_a) = \{u_iq^x|0 \le x \le n-1\}$ and hence
all the $\lambda_a$ belong to the same block. It remains to show
that if $\alpha, \beta$ are multipartitions and one of them is not
of the form $\lambda_a$, then they belong to distinct blocks.

\subsection{}
Now we suppose that we have two multipartitions $\alpha =
(\alpha^{(1)}, \ldots, \alpha^{(m)})$ and $\beta= (\beta^{(1)},
\ldots, \beta^{(m)})$ and $\cont(\alpha)=\cont(\beta)$. We will
show that if $k \neq i,j$ then $\alpha^{(k)}=\beta^{(k)}$.
\begin{lem}\label{lem4}
Let $k \neq i,j$. If $x \in \cont(\alpha^{(k)})$ then $x \in
\cont(\beta^{(k)})$.
\end{lem}
\begin{proof}
There is an integer $b$  with $-n+1 \le b \le n-1$ such that
$x=u_kq^b$. We consider the cases $b \ge 0$ and $b \le 0$
separately. In the case $b \ge 0$, we now prove by induction that
$x \notin \cont(\beta^{(\ell)})$ for any $\ell \neq k$. The proof
for $b \le 0$ is very similar, so we omit it.

For the base step, suppose $b=0$. Then $u_k \in
\cont(\alpha^{(k)})$. Hence $u_k$ is a residue of $\beta$. If $u_k
\in \cont(\beta^{(\ell)})$ where $\ell \neq k$ then $u_k=u_\ell
q^{c-r}$ for some column $c$ and row $r$ of $\beta^{(\ell)}$. But
clearly $-n<c-r<n$ which contradicts Lemma \ref{lem1}. Therefore
$u_k \notin \cup_{\ell \neq k}\beta^{(\ell)}$ and so $u_k \in
\cont(\beta^{(k)})$.

Now we do the inductive step. Suppose $b>0$. Suppose $u_kq^b$ is a
residue of $\beta^{(\ell)}$ with $\ell \neq k$. Then $u_kq^b=u_
\ell q^{c-r}$ for some $c,r$. So $u_k/u_\ell = q^{c-r-b}$. Since
$c-r<n$ and $b>0$, we have $c-r-b<n$. So by Lemma \ref{lem1},
$c-r-b \le -n$. Therefore, $r \ge n+c-b \ge n+1-b$. But
$\beta^{(\ell)}$ contains at least $r$ boxes, by definition of
$r$. So $|\beta^{(\ell)}| \ge n+1-b$.

Next, we note that since $u_kq^b$ is the residue of a node in
$\alpha^{(k)}$, this node must lie on the diagonal containing
$(b+1,1)$. So there are at least $b+1$ boxes in the first row of
$\alpha^{(k)}$ and hence there is a node in the first row of
$\alpha^{(k)}$ with residue $u_kq^{b-1}$. By induction on $b$,
this is also a residue of $\beta^{(k)}$. So there is a box in
column $b$ and row $1$ of $\beta^{(k)}$. Therefore, $|\beta^{(k)}|
\ge b$. So $|\beta| \ge |\beta^{(k)}|+|\beta^{(\ell)}| \ge n+1$, a
contradiction.
\end{proof}

\subsection{}
It follows from Lemma \ref{lem4} that if
$\cont(\alpha)=\cont(\beta)$ then
$\cont(\alpha^{(k)})=\cont(\beta^{(k)})$ for all $k \neq i,j$.
Then applying Lemma \ref{lem3}, we get $\alpha^{(k)}=\beta^{(k)}$.
It remains to deal with $\alpha^{(i)}$ and $\alpha^{(j)}$. The
proof of this case will be very similar to Lemma \ref{lem4}, but
slightly more complicated.
\subsection{}
Given multipartitions $\alpha = (\alpha^{(1)}, \ldots,
\alpha^{(m)})$ and $\beta= (\beta^{(1)}, \ldots, \beta^{(m)})$,
with $\cont(\alpha)=\cont(\beta)$, let $a_1$ be the length of the
first row of $\alpha^{(i)}$ and $a_2$ be the length of the first
column of $\alpha^{(j)}$ and define $b_1, b_2$ similarly for
$\beta$. First we prove a technical lemma.
\begin{lem}\label{lem5} Under our assumptions of Section \ref{stand},
suppose $a_1+a_2<n$. Then $u_iq^{a_1} \notin \cont(\alpha)$.
\end{lem}
\begin{proof}
First, we show that $u_iq^{a_1} \notin \cont(\alpha^{(k)})$ when
$k \neq i,j$. So let $k \neq i,j$ and suppose there is a node of
$\alpha^{(k)}$ with residue $u_iq^{a_1}$. Say this node lies in
column $c$ and row $r$ of $\alpha^{(k)}$. Then
$u_iq^{a_1}=u_kq^{c-r}$. So $u_i/u_k = q^{a_1-(c-r)}$. We show
that $a_1-(c-r)$ lies between $-n$ and $n$. If $a_1-(c-r) \ge n$
then $c+n \le r+a_1 \le n$, a contradiction. While if $a_1-(c-r)
\le -n$ then $c \ge n+a_1+r \ge n+1$, a contradiction. So
$-n<a_1-(c-r)<n$, which violates Lemma \ref{lem1}. Hence,
$u_iq^{a_1}$ is not a residue of $\alpha^{(k)}$.

Next, we show that $u_iq^{a_1}$ is not a residue of
$\alpha^{(i)}$. If it is, then there is a node in column $c$ and
row $r$ of $\alpha^{(i)}$ whose residue is
$u_iq^{a_1}=u_iq^{c-r}$. So $q^{a_1-(c-r)}=1$. So by Section
\ref{qorder}, if $a_1-(c-r) \neq 0$ then either $a_1-(c-r) \ge
2n-1$ or $a_1-(c-r)\le -(2n-1)$. If $a_1-(c-r) \le -(2n-1)$ then
$2n \le a_1+ r-1+2n \le c$, which is impossible. If $a_1-(c-r) \ge
2n-1$ then $c+2n \le a_1+r+1 \le n+2$, which is impossible if
$n>1$. Therefore, $a_1 =c-r$. But $c \le a_1$ and $r \ge 1$, so
this is also impossible. Therefore, $u_i q^{a_1}$ cannot be a
residue of $\alpha^{(i)}$.

The argument that $u_iq^{a_1}$ is not a residue of $\alpha^{(j)}$
is very similar. We use the fact that $a_1 < n-a_2$.
\end{proof}

%

\subsection{}
The claim of Section \ref{blocksclaim} follows from the next
lemma. We use the same notation as Section \ref{lem5}.
\begin{lem}\label{lem6}
Under the assumptions of Section \ref{stand}, if $a_1+a_2<n$ then
if $x \in \cont(\alpha^{(i)})$ then $x \in \cont(\beta^{(i)})$.
\end{lem}
\begin{proof}
By Lemma \ref{lem4}, $\cont(\alpha^{(k)})=\cont(\beta^{(k)})$ for
$k \neq i,j$. Therefore, by Lemma \ref{lem2}, we get
$\cont(\alpha^{(i)}) \cup \cont(\alpha^{(j)}) =
\cont(\beta^{(i)})\cup\cont(\beta^{(j)})$. This is a disjoint
union, so it suffices to show that if $x \in \cont(\alpha^{(i)})$
then $x \notin \cont(\beta^{(j)})$.

If $x \in \cont(\alpha^{(i)})$ then $x =u_iq^b$ for some $b$ with
$-n+1 \le b \le n-1$. As in the proof of Lemma \ref{lem4}, we
consider the cases $b \ge 0$ and $b \le 0$ separately. We give the
proof only for the $b \ge 0$ case. The proof is by induction on $b$.

For the base step, if $b=0$ then $u_i$ is a residue of
$\alpha^{(i)}$. If this is a residue of $\beta^{(j)}$, then it has
the form $u_i=u_iq^{n-1}q^{c-r}$ for some $c,r$. So
$q^{n-1+c-r}=1$. Now, $n-1+c-r \ge 0$. If $n-1+c-r \ge 2n-1$ then
$c-r \ge n$ which is impossible. So $n-1+c-r=0$. Hence, $c=1,
r=n$, and $\beta^{(j)}$ must be a column of $n$ boxes. But then
$\cont(\beta^{(j)})=\{u_iq^{n-1},u_iq^{n-1},\ldots,u_iq,u_i\}$.
Since $0 \le a_1 < n$, we have $u_iq^{a_1} \in
\cont(\beta^{(j)})=\cont(\beta)=\cont(\alpha)$, which contradicts
Lemma \ref{lem5}. Therefore $u_i$ must be a residue of
$\beta^{(i)}$, which proves the base step.

For the inductive step, suppose $b>0$ and $u_iq^b$ is a residue of
$\alpha^{(i)}$. If $u_iq^b$ is a residue of a node in column $c$ and
row $r$ of $\beta^{(j)}$, then $u_iq^b=u_iq^{n-1}q^{c-r}$. So
$q^{c-r+n-1-b}=1$. Since $c-r<n$ and $b>0$, we have $c-r-b<n$. So
$c-r-b+n-1<2n-1$. Therefore, either $c-r-b+n-1=0$ or $c-r-b+n-1 \le
-(2n-1)$. If the latter holds then $c+3n \le r+b+2 \le 2n+1$ since
we may take $b \le n-1$. Hence $1+n \le c+n \le 1$, a contradiction.
We therefore get $c-r-b+n-1=0$. So $r \ge n-b$. But $\beta^{(j)}$
has at least $r$ nodes. Therefore, $|\beta^{(j)}| \ge n-b$ and has
at least $n-b$ rows. But since $u_iq^b \in \cont(\alpha^{(i)})$, we
get $u_iq^{b-1} \in \cont(\alpha^{(i)})$, as in the proof of Lemma
\ref{lem4}. By induction on $b$, $u_iq^{b-1} \in
\cont(\beta^{(i)})$. So, as in the proof of Lemma \ref{lem4}, there
is a box in row $1$ and column $b$ of $\beta^{(i)}$. Therefore,
$|\beta^{(i)}| \ge b$ and $\beta^{(i)}$ has at least $b$ columns. So
$\beta =\lambda_b$ in the notation of Section \ref{blocksclaim}.
Therefore $\cont(\beta)=\{u_i, qu_i, \ldots, q^{n-1}u_i\}$. So
$u_iq^{a_1} \in \cont(\beta)= \cont(\alpha)$. This contradicts Lemma
\ref{lem5}. Therefore, $u_iq^b$ must be a residue of $\beta^{(i)}$
and this proves the inductive step.
\end{proof}

\subsection{}
Now suppose we have a multipartition $\alpha$ not of the form
$\lambda_a$. Suppose $\beta \neq \alpha$. We show that
$\cont(\alpha) \neq \cont(\beta)$. Indeed, if $\beta \neq
\lambda_b$ for any $b$, then by Lemmas \ref{lem4} and \ref{lem6},
$\cont(\alpha^{(k)})=\cont(\beta^{(k)})$ for all $k$. Therefore,
by Lemma \ref{lem3}, $\alpha^{(k)}=\beta^{(k)}$ for all $k$, so
$\alpha=\beta$, a contradiction. On the other hand, if $\beta =
\lambda_b$ for some $b$, then $u_i q^{a_1} \in \cont(\beta)
\setminus \cont(\alpha)$ by Lemma \ref{lem5}. So $\cont(\alpha)
\neq \cont(\beta)$.

 Therefore, $S^\alpha$ is the unqiue Specht module
in its block. Furthermore, $\{S^{\lambda_a}|0 \le a \le n\}$ form
a block, by the same reasoning.

\subsection{}\label{parts3and4}
We get that there is one block of the Hecke algebra containing
$n+1$ of the Specht modules, and all the other blocks are
singletons. Hence, there are $|\textsf{Irrep}(W)|-n$ blocks. By
\cite[Corollary 5.18]{GGOR}, the blocks of $\cO$ are in bijection
with blocks of $\mathcal{H}$ and hence $\cO$ also has
$|\textsf{Irrep}(W)|-n$ blocks. We work in the category
$\cO(H_{\kappa'})$. Now by \cite[Theorem 2.3]{CE}, there is a
BGG-resolution of $\tilde{Y_c}$, ie. an exact sequence
$$0 \leftarrow \tilde{Y_c} \leftarrow M(\mathsf{triv}) \leftarrow
M(\h_q) \leftarrow \cdots \leftarrow M(\wedge^n\h_q) \leftarrow
0$$ As the classes $[M(\tau)]$ form a basis of the Grothendieck
group $K_0(\cO)$, none of the maps in this sequence can be zero,
and hence all the $L(\wedge^i\h_q)$ belong to the same block.
There are $n+1$ simples in this block and hence by counting we see
that all the other blocks must be singletons. Using the fact that
simple objects in $\cO$ have no self-extensions (\cite[Proposition
1.12]{BEG2}), we get that these blocks are semisimple. Translating
back to category $\cO(H_\kappa)$, we get parts $(3)$ and $(4)$ of
Theorem \ref{mainthm}.

\subsection{}\label{compo}
It remains to compute the composition multiplicities in the one
nontrivial block $\cO^\wedge$. Again we work in the category
$\cO(H_{\kappa'})$. \cite[Proposition 5.21(ii)]{GGOR} tells us
that each $L(\wedge^i\h_q)$, $i>0$ is a submodule of a standard
module. Write $L_i = L(\wedge^i\h_q)$ and $M_i = M(\wedge^i\h_q)$.
Let $R_i$ be the radical of $M_i$. We cannot have a nonzero map
$L_i \rightarrow M_j$ if $j
>i$ by \cite[Section 2.5 (32)]{DO} (briefly, this is because a
calculation very similar to \cite[Lemma 8.3]{gmpn} shows that the
number denoted $c_{\wedge^t\h_q}(k)$ in \cite{DO} equals $-tN$ for
some $N \in \mathbb{N}$ which is independent of $t$. Thus if
$[M_j:L_i] \neq 0$ then $-jN+iN \in \mathbb{N}$ and so $i>j$) and
so $L_1$ is a submodule either of $M_0$ or $M_1$. It can't be a
submodule of $M_1$ because $[M_1:L_1]=1$ so $L_1 \hookrightarrow
M_0$. So $L_1 \hookrightarrow R_0$. But $R_0$ is a quotient of
$M_1$ since it is the image of $M_1 \rightarrow M_0$ (this follows
from the fact that $\tilde{Y}_c \cong L_0$, proved in
\cite[Section 8.2]{gmpn}), hence $[R_0:L_1]=1$. If we had
$[R_0:L_i] \neq 0$ for some $i
> 1$ then $R_0$ would have $L_i$ as a quotient for some $i>1$.
Therefore, so would $M_1$. But $M_1$ has a unique simple quotient
$L_1$. Therefore, it is impossible to have $[R_0:L_i] \neq 0$ for
$i >1$ and we conclude that $R_0=L_1$.

\subsection{}\label{part5}
We have shown that the composition factors of $M_0$ are $L_0$ and
$L_1$. To conclude the argument, we show by induction that the
composition factors of $M_i$ are $L_i$ and $L_{i+1}$. Consider
first $L_{i+1}$. Then $L_{i+1}$ is a submodule of some $M_j$. We
cannot have $j \ge i+1$, and by induction, we cannot have $j<i$.
Hence, $L_{i+1}$ is a submodule of $M_i$ and so $L_{i+1}
\hookrightarrow R_i$. Now $R_i = \ker(M_i \rightarrow M_{i-1})$ by
induction and so $R_i$ is a quotient of $M_{i+1}$. Therefore,
$[R_i:L_{i+1}]=1$. If there was a $j>i+1$ with $[R_i:L_j] \neq 0$
then we would have that for some $j>i+1$, $L_j$ would be a
quotient of $R_i$ and hence a quotient of $M_{i+1}$, contradicting
the fact that $M_{i+1}$ has a unique simple quotient. Therefore,
$R_i =L_{i+1}$ and we are done. This proves part $(5)$ of Theorem
\ref{mainthm}.
\section{Characterisations of separating simples}\label{parameterchoices}
Now that we have completed the proof of Theorem \ref{mainthm}, let
us turn our attention to the question of when $\KZ$ separates
simples.
\begin{thm}\label{ssintermsofAK}
The following are equivalent
\begin{enumerate}
\item $\KZ$ separates simples. \item If $q$, $u_1, \ldots, u_m$
are the parameters of the Ariki-Koike algebra $\calH$, then
$(q+1)\prod_{i<j} (u_i -u_j) \neq 0$, and furthermore,
$$\#\{ \tau \in \irrep(W) : \ L(\tau)|_{\hreg} \neq 0 \} \ge n-1.$$
\item The algebra $\calH$ has at least $|\irrep(W)|-1$
nonisomorphic simple modules.
\end{enumerate}
\end{thm}
\begin{proof}
First, we show that (2) implies (1).
We must show that if
$L(\sigma)|_{\mathfrak{h}^\mathrm{reg}} \cong
L(\tau)|_{\mathfrak{h}^\mathrm{reg}} \neq 0$ then $\sigma=\tau$.
Suppose then that $L(\sigma)|_{\mathfrak{h}^\mathrm{reg}} \cong
L(\tau)|_{\mathfrak{h}^\mathrm{reg}} \neq 0$. By \cite[Proposition
5.21(ii)]{GGOR}, there exists a standard module $M(\lambda)$ such
that $L(\sigma) \hookrightarrow M(\lambda)$. Let $t =
\dim(\Hom(L(\sigma),M(\lambda)))$. Then $M(\lambda)$ must
have $t$ submodules isomorphic to $L(\sigma)$, because the
only automorphisms of $L(\sigma)$ are the scalars. 
Therefore, 
$L(\sigma)^{\oplus t} \subset M(\lambda)$ and $M(\lambda)$ has no
submodule isomorphic to $L(\sigma)^{\oplus (t+1)}$. Now since
$L(\sigma)|_{\mathfrak{h}^\mathrm{reg}} \cong
L(\tau)|_{\mathfrak{h}^\mathrm{reg}}$, we have
$\Hom(L(\tau)|_{\mathfrak{h}^\mathrm{reg}},
M(\lambda)|_{\mathfrak{h}^\mathrm{reg}}) =
\Hom(L(\sigma)|_{\mathfrak{h}^\mathrm{reg}},
M(\lambda)|_{\mathfrak{h}^\mathrm{reg}}) \neq 0$ and hence by
\cite[Proposition 5.9]{GGOR}, $\Hom(L(\tau), M(\lambda)) \neq 0$
(using the condition on the parameters). Therefore, $M(\lambda)$
has a submodule isomorphic to $L(\tau)$ and hence a submodule
isomorphic to $L(\tau)+L(\sigma)^{\oplus t}$. This sum must be
direct if $L(\sigma) \ncong L(\tau)$, hence $M(\lambda)$ has a
submodule $L(\tau) \oplus L(\sigma)^{\oplus t}$ and
$M(\lambda)|_{\mathfrak{h}^\mathrm{reg}}$ has a submodule
$L(\tau)|_{\mathfrak{h}^\mathrm{reg}} \oplus
L(\sigma)|_{\mathfrak{h}^\mathrm{reg}}^{\oplus t} =
L(\sigma)|_{\mathfrak{h}^\mathrm{reg}}^{\oplus (t+1)}$. Therefore,
$$\mathrm{dim}(\Hom(L(\sigma)|_{\mathfrak{h}^\mathrm{reg}},
M(\lambda)|_{\mathfrak{h}^\mathrm{reg}})) \ge t+1$$ and so
$\mathrm{dim}(\Hom(L(\sigma), M(\lambda))) \ge t+1$, a
contradiction. So $L(\sigma) \cong L(\tau)$ and hence
$\sigma=\tau$.

Next, (1) implies (3) by Section \ref{beginning}.

Finally, to show (3) implies (2), note that under the hypothesis
that $\calH$ has $|\irrep(W)|-1$ simple modules, it has already
been shown 
that $[n]_q! \neq 0$, hence $q \neq -1$ since we assume $n \ge 2$, and 
that $u_i \neq u_j$ for all $i
\neq j$, so the condition on the parameters holds. Furthermore, by
the essential surjectivity of $\KZ$, if $\calH$ has
$|\irrep(W)|-1$ simple modules then, because $\KZ$ is essentially
surjective on objects and exact, there are at least
$|\irrep(W)|-1$ of the $L(\tau)$ with $\KZ(L(\tau)) \neq 0$ and
hence with $L(\tau)|_{\hreg} \neq 0$.
\end{proof}
\section{The Ariki-Koike algebra in the almost-semisimple case}
\subsection{}
We close this section by using the facts proved about category
$\cO$ in Theorem \ref{mainthm} to prove a theorem about the Hecke
algebra which does not mention the Cherednik algebra in its
hypothesis or conclusion. This theorem is an example of a general
philosophy suggested by Rouquier in \cite{rouquiersurvey} of using
the Cherednik algebra and the $\KZ$ functor as a tool to prove
theorems about Hecke algebras.

It is well-known that $\calH_\kappa$ is semisimple if and only if the
number of irreducible modules $|\irrep(\calH_\kappa)|$ of
$\calH_\kappa$ equals the number of irreducible modules of $\C W$,
and that in this case $\calH_\kappa \cong \C W$. So the property
of having $|\irrep(W)|$ simple modules determines the algebra
$\calH_\kappa$ up to isomorphism. We show that the property of
having $|\irrep(W)|-1$ simple modules also determines
$\calH_\kappa$ up to isomorphism.
\begin{thm}\label{akassc}
Suppose $\calH_\kappa$ and $\calH_\mu$ are Ariki-Koike algebras
corresponding to some parameters $\kappa , \mu \in \C^m$ and
that $|\irrep (\calH_\kappa)|=|\irrep (\calH_\mu)|= |\irrep
(W)|-1$. Then there is an isomorphism of algebras $\calH_\kappa
\cong \calH_\mu$.
\end{thm}
\begin{proof}
By \cite[Theorem 5.15]{GGOR}, there is an algebra isomorphism $\calH_\kappa \cong \mathrm{End}_{\cO}(P_{\KZ})^{opp}$ where $$P_{\KZ} = \bigoplus_{\tau \in \irrep(W)} \dim \KZ (L(\tau)) P(\tau).$$ Here, $P(\tau)$ is the projective cover of $L(\tau)$. The strategy of the proof is to calculate $P_\KZ$ in the case where $\KZ_\kappa$ separates simples, and show that its endomorphism ring can be written in a way that does not depend on $\kappa$. We work in the category $\cO = \cO_\kappa$ and write $\KZ=\KZ_\kappa$, $M(\tau) = M_\kappa (\tau)$, and so forth. By Theorem \ref{mainthm}, there is a linear representation $\chi$ of $W$ with $\cO=\cO^\wedge \oplus \cO^{ss}$, where $\cO^\wedge$ is the subcategory of $\cO$ generated by $\{ L(\wedge^i \h_q \otimes \chi) : 0 \le i \le n\}$. Let $\lambda^i = \wedge^i\h_q \otimes \chi$ and let $S= \{ \lambda^i : 0 \le i \le n\}$. Write $M_i =M(\lambda^i)$, $L_i =L(\lambda^i)$ and $P_i = P(\lambda^i)$ (the projective cover of $L_i$). 

For $\sigma, \tau \in \irrep(W)$, we have in general 
\begin{align*}
\dim \mathrm{Hom} (P(\sigma), P(\tau)) &= [P(\tau):L(\sigma)] \\
&= \sum_\gamma [P(\tau):M(\gamma)][M(\gamma):L(\sigma)]\\
&= \sum_\gamma [M(\gamma):L(\tau)][M(\gamma):L(\sigma)]\\
&= \sum_{\gamma \in S} [M(\gamma):L(\tau)][M(\gamma):L(\sigma)]
+ \sum_{\gamma \notin S}[M(\gamma):L(\tau)][M(\gamma):L(\sigma)]
\end{align*}
If $\gamma \notin S$ then $M(\gamma)=L(\gamma)$, so we get
$$\dim \mathrm{Hom}(P(\sigma),P(\tau)) = \sum_{i=0}^n[M_i:L(\tau)][M_i:L(\sigma)]+\sum_{\gamma \notin S}\delta_{\gamma \tau} \delta_{\gamma \sigma}.$$
Now, if $\sigma \notin S$ or $\tau \notin S$, this sum must be $\delta_{\sigma \tau}$. Otherwise, $\sigma, \tau \in S$ and so $\sigma = \lambda^a$, $\tau = \lambda^b$ for some $a,b$. We get
$$\dim \mathrm{Hom}(P(\lambda^a),P(\lambda^b)) = \sum_{i=0}^n[M_i:L_a][M_i:L_b]$$
which equals $2$ if $a=b$ and $1$ if $|a-b|=1$ and $0$ otherwise. So we get
$$\dim \mathrm{Hom}(P(\sigma),P(\tau)) =
\begin{cases}
2 & \text{if } \sigma = \tau \in S\\
1 & \text{if } \sigma = \tau \notin S\\
1 & \text{if } \{ \sigma, \tau \}=\{ \lambda^a,\lambda^{a+1} \}, 0 \le a \le n-1\\
0 & \text{otherwise}\end{cases}$$
The ring $\mathrm{End}_\cO(P_\KZ)$ is a matrix algebra with entries in the various Hom-spaces $\mathrm{Hom}(P(\sigma), P(\tau))$. We calculate the  multiplication relations between basis elements of the $\mathrm{Hom}(P(\sigma), P(\tau))$ and show that these relations do not depend on $\kappa$. It will follow that the structure constants of  $\mathrm{End}_\cO(P_\KZ)$ do not depend on $\kappa$, which will prove the theorem provided that the multiplicaity of each $P(\tau)$ in $P_{\KZ}$ is also independent of $\kappa$. But in our situation $P_\KZ = \oplus_{\tau \notin S} (\dim \tau)\cdot P(\tau) \oplus \left( \oplus_{1 \le i \le n} {n-1 \choose i-1} P_i \right)$ since $\dim \KZ (L_i)= {n-1 \choose i-1}$, as can be readily shown using induction on the BGG-resolution of $L_0$ and the fact that $\dim \KZ (M(\tau)) = \dim (\tau)$ for all $\tau$.

By BGG reciprocity, we have $[P_i:M_i]=[P_i:M_{i-1}]=1$, and $[P_i:M(\sigma)]=[M(\sigma):L_i]=0$ if $\sigma \neq \lambda^i, \lambda^{i-1}$. Therefore, the factors in any filtration of $P_i$ by standard modules are $M_i$ and $M_{i-1}$. But by \cite[Corollary 2.10]{GGOR}, $P_i$ has a filtration by standard modules with $M_i$ as the top factor, so $P_i$ may be described as $P_i=\begin{smallmatrix}M_i \\ M_{i-1}\end{smallmatrix}$, meaning that there is a series $0=P_i^0 \subset P_i^1 \subset P_i^2 =P_i$ with $P_i^1 \cong M_{i-1}$ and $P_i^2/P_i^1 \cong M_i$. We may write the resulting composition series of $P_i$ as
$$P_i= \begin{matrix} L_i \\ L_{i+1} \\ L_{i-1} \\ L_i \end{matrix}$$ 
This description of $P_i$ makes it easy to write down the nontrivial maps $P_i \rightarrow P_{i}$. 

First, there are two obvious maps $P_i \rightarrow P_i$, namely the identity map $\mathrm{id}_i$ and the map $\xi_i$ which is projection onto the top composition factor $L_i$ followed by inclusion. Note that $\xi_i^2=0$ and therefore $\mathrm{End}_\cO(P_i) = \C[\xi_i]/(\xi_i^2)$, since we have already shown that $\dim \mathrm{Hom}(P_i,P_i)=2$.

Next, we describe the map $P_i \rightarrow P_{i+1}$. This is a map $\begin{smallmatrix} M_i \\ M_{i-1} \end{smallmatrix} \rightarrow \begin{smallmatrix} M_{i+1} \\ M_i \end{smallmatrix}$. So we may construct a map $f_{i,i+1}: P_i \rightarrow P_{i+1}$ by factoring out the copy of $M_{i-1}$ and then embedding $M_i$ in $P_{i+1}$. This map is nonzero, so $\mathrm{Hom}(P_i, P_{i+1}) = \C f_{i,i+1}$, $1 \le i \le n-1$. 

Now we describe the map $P_i \rightarrow P_{i-1}$, $n \ge i \ge 2$.
By \cite[Proposition 5.2.1 (ii)]{GGOR}, $P_i \supset L_i$ is injective and therefore $P_i$ contains the injective envelope $I_i = I(\lambda^i)$ of $L_i$. 
Therefore, since $P_i$ is indecomposable,  $P_i =I_i$. Now, category $\cO$ contains a \emph{costandard module} $\nabla(\tau) \supset L(\tau)$ for every $\tau \in \irrep(W)$, with $[\nabla(\tau)]=[M(\tau)]$ in $K_0(\cO)$. Write $\nabla_i = \nabla(\lambda^i)$. Then $L_i \subset \nabla_i$, so $\nabla_i$ has a composition series of the form $\nabla_i = \begin{smallmatrix} L_{i+1} \\ L_i \end{smallmatrix}$. Furthermore, $\nabla_i \subset I_i$ and so $I_i$ has a filtration by costandard modules of the form $I_i = \begin{smallmatrix} \nabla_{i-1} \\ \nabla_i \end{smallmatrix}$ (the existence of such a filtration follows from \cite[Definition 3.1, Axiom (c)]{cps}). Since $I_i = P_i$, to get a map $\begin{smallmatrix} \nabla_{i-1} \\ \nabla_i \end{smallmatrix}=P_i \rightarrow P_{i-1}=\begin{smallmatrix} \nabla_{i-2} \\ \nabla_{i-1} \end{smallmatrix}$, we may factor out the copy of $\nabla_i$ and then embed $\nabla_{i-1}$ in $P_{i-1}$. This gives a nonzero map $f_{i,i-1}$, and therefore $\mathrm{Hom}(P_i,P_{i-1})=\C f_{i,i-1}$.
In particular, this shows that the image of $f_{i,i-1}$ has length $2$. 

Now we calculate multiplication relations between the various $f_{i,i+1}$, $f_{i,i-1}$ and $\xi_i$. First, it is immediate from the definitions that $\xi_{i+1} f_{i,i+1} = f_{i,i+1} \xi_i=0$. We need to do a little more work to show that the same holds for $f_{i,i-1}$. Take the description of $I_i$ as $I_i = \begin{smallmatrix} \nabla_{i-1} \\ \nabla_i \end{smallmatrix}$. Then $I_i$ has a composition series 
$$I_i = \begin{matrix} L_i \\ L_{i-1} \\ L_{i+1} \\ L_i \end{matrix}$$
So there is a map $\zeta_i : I_i \rightarrow I_i$ defined by projection onto the top composition factor $L_i$ followed by the embedding $L_i \hookrightarrow I_i$. Clearly, $\zeta_i f_{i-1,i} = f_{i-1,i} \zeta_{i-1}=0$. But since $P_i=I_i$, we may regard $\zeta_i$ as a map $P_i \rightarrow P_i$. Therefore, there are $a,b \in \C$ with $\zeta_i = a \mathrm{id}_i + b\xi_i$. Since $\zeta_i^2=0$, we get $a^2=0$ and hence $\zeta_i$ is a nonzero multiple of $\xi_i$. This shows that $\xi_i f_{i-1,i}=f_{i-1,i}\xi_{i-1}=0$.

Finally, we need to calculate $f_{i+1,i}f_{i,i+1}$ and $f_{i-1,i}f_{i,i-1}$. Consider first $f_{i-1,i}f_{i,i-1}$. By the 
definition of $f_{i,i-1}$ above, we have $[\mathrm{im}(f_{i,i-1}):L_i]\neq 0$. Hence, $\mathrm{im}(f_{i,i-1})$ cannot be contained in the submodule of $P_{i-1}$ isomorphic to $M_{i-2}$, and therefore $f_{i-1,i}f_{i,i-1}$ must be nonzero. Since $f_{i-1,i}f_{i,i-1} \xi_i=0$, $f_{i-1,i}f_{i,i-1}$ must be a nonzero multiple of $\xi_i$. Let us replace $\xi_i$ by $f_{i-1,i}f_{i,i-1}$. So we may assume that $f_{i-1,i}f_{i,i-1} = \xi_i$, and this does not change any of the relations which have already been calculated. Now consider $f_{i+1,i}f_{i,i+1}$. We show that this composition is nonzero. Indeed, the image $\mathrm{im}(f_{i,i+1})$ has composition factors $L_i$ and $L_{i+1}$. If $f_{i+1,i}f_{i,i+1}$ were zero, then we would get that $\mathrm{im}(f_{i+1,i})$ could only have composition factors $L_{i+1}$ and $L_{i+2}$. But we have shown that $\mathrm{im}(f_{i+1,i})$ has length $2$, and $[P_i:L_{i+2}]=0$, a contradiction. Therefore, $f_{i+1,i}f_{i,i+1} \neq 0$ and so there is a nonzero $b_{i,i+1} \in \C$, $n-1 \ge i \ge 1$, such that 
$$f_{i+1,i}f_{i,i+1} = b_{i,i+1} \xi_i = b_{i,i+1} f_{i-1,i}f_{i,i-1}.$$
It remains to do some rescaling. Let 
\begin{align*}
\xi_i' &= \frac{1}{b_{12}b_{23} \cdots b_{i-1,i}} \xi_i, \quad &1 \le i \le n \\
f_{i,i-1}' &= f_{i,i-1} \quad &2 \le i \le n \\
f_{i,i+1}' &= \frac{1}{b_{12}b_{23} \cdots b_{i,i+1}}f_{i,i+1} \quad &1 \le i \le n-1.
\end{align*}
Then we have the following relations:
\begin{align}\label{frels}\notag
\xi_i' f_{i-1,i}' &= f_{i-1,i}' \xi_{i-1}' =0\\ \notag
\xi_{i+1}' f_{i,i+1}' &= f_{i,i+1}' \xi_i' =0 \\ 
f_{i-1,i}'f_{i,i-1}' &= f_{i+1,i}'f_{i,i+1}' =\xi_i'. 
\end{align}
These are the only nontrivial relations between the various $\mathrm{Hom}(P(\sigma),P(\tau))$. This shows that we may choose a basis of $\mathrm{Hom}(P(\sigma),P(\tau))$ for each $\sigma, \tau$ such that the composition relations between the basis elements are independent of $\kappa$. Hence, we may choose a basis of the algebra $\mathrm{End}_\cO(P_\KZ)$ such that the structure constants are independent of $\kappa$. This proves the theorem.
\end{proof}
\begin{rem}
By variations on the arguments given in the above proof, it is possible to show that 
$$\ext^1(L_i,L_j)= \begin{cases} 1 & j = i+1, i-1 \\
0 & \text{otherwise} \end{cases}$$
and so the composition series of $P_i$ may be written more symmetrically as $P_i=I_i= \begin{smallmatrix} L_i \\ L_{i-1} \oplus L_{i+1} \\ L_i \end{smallmatrix}$.
\end{rem}
Note that since we have shown earlier that the Ariki-Koike algebra has $|\irrep(W)|-n$ blocks, by counting we get that the algebra $B_n := \mathrm{End}_\cO (\oplus_{i=1}^n {n-1 \choose i-1} P_i)$ is a block of the Ariki-Koike algebra. From the relations (\ref{frels}), it is clear that $B_n$ is independent both of $\kappa$ and $m$. Furthermore, we may extend this description of the unique non-semisimple block to $m=1$. This is because in the $m=1$ case, the Cherednik algebra only depends on one parameter $\kappa_{00}$ (denoted $c$ in \cite{BEG2}). We write the Hecke algebra as $\calH_c (S_n)$, with parameter $q= e^{2 \pi i c}$. The simple modules of $\calH_c(S_n)$ are in bijection with $e$--restricted partitions $\lambda$ of $n$, where $e$ is the muliplicative order of $q$ in $\C^*$, and a partition $\lambda$ is said to be $e$--restricted if $\lambda_i -\lambda_{i+1}<e$ for all $i \ge 1$. It is clear from this decription that $\calH_c$ has $|\irrep(S_n)|-1$ simple modules if and only if $e=n$ if and only if $c = \frac{r}{n}$ with $(r,n)=1$. In this case, Theorem \ref{mainthm} holds without change by various results of \cite[Section 3]{BEG2}, and the proof of Theorem \ref{akassc} also goes through without change in the case $m=1$. So we have the following corollary.
\begin{cor}
Let $\ell_1, \ell_2 \ge 1$ and for $i=1,2$ let $\kappa_i \in \C^{\ell_i}$ and suppose $\calH_{\kappa_i}(G(\ell_i,1,n))$ has $|\irrep(G(\ell_i,1,n))|-1$ simple modules. Then the unique nonsemisimple blocks of $\calH_{\kappa_1}(G(\ell_1,1,n))$ and $\calH_{\kappa_2}(G(\ell_2,1,n))$ are isomorphic algebras. In particular, they are isomorphic to the principal block $B_n$ of $\calH_{1/n}(S_n)$.
\end{cor}
\begin{rem}
The representation theory of the algebra $B_n$ is described in \cite[5.3]{BEG2} and \cite[3.2]{EN}.
\end{rem}

\bibliographystyle{alpha}
\bibliography{penguin}
\end{document}